\documentclass[twoside]{article}

\usepackage{PRIMEarxiv}

\usepackage{indentfirst}
\usepackage[colorlinks,
linkcolor=blue,
anchorcolor=black,
urlcolor=blue,
citecolor=red]{hyperref}

\usepackage{url}

\usepackage{amssymb}
\usepackage{amsthm}
\usepackage[]{amsmath}
\usepackage{cuted}
\allowdisplaybreaks

\usepackage{enumerate}

\usepackage{mathptmx}

\usepackage{subfigure}
\usepackage{graphicx}
\usepackage{epsfig}

\usepackage{multirow}
\usepackage{array}
\usepackage{floatrow}
\usepackage{makecell}

\usepackage{booktabs}
\usepackage[font=small,labelfont=bf,singlelinecheck=false]{caption}
\usepackage[justification=raggedright]{caption}

\usepackage{listings}

\usepackage{lipsum}

\usepackage[utf8]{inputenc} % allow utf-8 input
\usepackage[T1]{fontenc}    % use 8-bit T1 fonts
\usepackage{hyperref}       % hyperlinks
\usepackage{url}            % simple URL typesetting
\usepackage{booktabs}       % professional-quality tables
\usepackage{amsfonts}       % blackboard math symbols
\usepackage{nicefrac}       % compact symbols for 1/2, etc.
\usepackage{microtype}      % microtypography
\usepackage{lipsum}
\usepackage{fancyhdr}       % header
\usepackage{graphicx}       % graphics
\graphicspath{{media/}}     % organize your images and other figures under media/ folder

\title{Exploring incentive strategies and predicting development trends for new energy vehicles}

\author{
	 Tao Jin \\
	 Key Laboratory of Electronic Information of State Ethnic Affairs Commission \\
	College of Electrical Engineering\\
	Southwest Minzu University \\
	Chengdu, Sichuan, China\\
	\texttt{jintao0612@stu.swun.edu.cn} \\
	\And
   	Yulian Jiang, Xingwen Liu \\
 	College of Electrical Engineering\\
 	Southwest Minzu University \\
	Chengdu, Sichuan, China\\
}

\begin{document}
\maketitle

\begin{abstract}
To facilitate new energy vehicles (NEVs), we construct a game model between vehicle manufacturers and consumers to explore their interactions. In the model, we propose the Expectation Supply-Demand Game (ESDG), construct the consumer purchasing decision-making process with feedback and analyse the stability of the system under different feedback factors. We processes the data of the model in numerical simulation through Min-Max normalisation and predicts the development of NEVs. The results show that: (1) An evolutionary stabilisation strategy (ESS) emerges in the evolutionary game model with the introduction of feedback. (2) The Min-Max normalisation method is conducive to the accuracy of the model. (3) Excessive advertising and marketing may cause consumer boredom. (4) The establishment of an appropriate battery compensation and replacement insurance is conducive to the development of NEVs. (5) The production and sales ratio of China's NEVs is predicted to reach 37.2\% and 36.9\% respectively in 2024. 
\end{abstract}

% keywords can be removed

\section{Introduction}
	As China's environmental and energy problems become more and more prominent, NEVs have attracted widespread attention as a sustainable green energy industry, and they play an important role in China's energy structure transformation. With the urgency of developing the NEV industry, although the Chinese government, enterprises and other relevant organisations have implemented various initiatives to further promote NEVs, the current situation of production and sales of NEVs still has a lot of room for growth. It is necessary to continue to implement more measures to stimulate the input and support of the relevant parties.

	The development of the NEV industry can	be credited to consumers and vehicle manufacturers. The direct stakeholders in the development of the NEV industry are consumers and vehicle manufacturers. Analyzing their interactions can provide insights into the future of NEVs. And evolutionary game theory is a useful tool for exploring economic interactions among agents in a dynamic, finite and rational two-player game \cite{1998On,Evolutionary}. Thus, this study aims to conduct an evolutionary game model to simulate the dynamic interactions between these parties, which can help explore strategies that favour the promotion of NEVs. To be specific, we proposed NEV battery payout and replacement insurance measures and explored the way of implementing external demand stimulus. Secondly, we propose a base model called the ESDG to understand the behavioural interactions between vehicle manufacturers and consumers. Meanwhile, we introduced interdisciplinary mathematical methods to explore their impact on improving evolutionary model accuracy, including consumer decision-making process with feedback and Min-Max normalisation. These methods have previously been widely used in the fields of psychology and machine learning, which means that this study can go some way to filling some of the research gaps, as will be described in Section \ref{section2}. Finally, We will analyse the equilibrium stability of the constructed evolutionary game model, which will also be verified by a sensitivity analysis of the numerical simulation. Our model achieved the reproduction and prediction of the development path of NEVs at the end. 

Through this study, the new fundamental game model we proposed can enrich game theory with a new theoretical framework to try to understand a specific real-world situation. The application areas of some mathematical tools can also be extended in this study, which provides different ideas and methods for the construction of mathematical models. In addition, this study can provide certain scientific references for relevant subjects to support NEVs. As such, our research aims to explore the following specific questions:  
\begin{enumerate}[\quad (1)]	
	\item How can the relationship between vehicle manufacturers and consumers be represented by a basic game model?
	\item How does the introduction of Min-Max normalisation affect the numerical simulation of evolutionary games?
	\item How does the consumer decision-making process with the introduction of a feedback factor affect the stability of evolutionary game systems, and what strategic insights can be drawn from the process?
	\item How does the consumer decision-making process with the introduction of a feedback factor affect the stability of evolutionary game systems, and what strategic insights can be drawn from the process?
	\item Are battery payout and replacement insurance measures beneficial to the evolutionary path of NEVs?
\end{enumerate}

\section{Related literature}
\label{section2}

As a clean energy industry, the development of NEVs has become an important research topic for scholars. Some scholars have discussed the current situation in relation to government and business, such as subsidies \cite{WOS:000793460600001}, industrial policies \cite{WOS:000609019000020}, the mapping relationship among strategies \cite{Bai2017Incentive}, and the effectiveness of related strategies \cite{WOS:000927730400001,WOS:000862763300005,JIN2023127677}. Other studies have primarily explored financial subsidies, with the latest studies having attempted to introduce specific proposals, including bans on petrol cars \cite{WOS:000810096300005}, carbon trading \cite{WOS:000832135200001,WOS:000848270400001}, battery recycling \cite{WOS:000930068000001}, etc. Yet, none of them has proposed strategy improvements to explore specific implementation options. Different strategies deserve to be explored. And the main subjects of the game involved in these studies are the government and the vehicle manufacturers. And consumers, as a key factor in the development of NEVs, it is necessary for them to be considered as a subject in the game model. Based on this, we tried to construct an evolutionary game model that includes consumers. Meanwhile, the impact of the government's initial strategy on consumer behaviour is relatively small, and with the added dimension of the number of participants in the three-way evolutionary game model, the complexity of the model increases, leading to certain limitations in its deployment and operation in a computing environment. To pursue a simpler and more accurate mathematical model structure, we mainly focused on the game interaction between vehicle manufacturers and consumers to propose different incentives for the development of NEVs. In addition, the research related to NEVs also has included many evolutionary game models. For example, \cite{2018Evolutionary} simulated the behavioral changes of the government and enterprises when the subsidies were withdrawn. \cite{2018Critical} studied the dynamic process of ``deception regulation" between the government and vehicle manufacturers. The consideration of the consumer in these studies has been directed towards a particular stage rather than a holistic process analysis, which to some extent is not conducive to the construction of mathematical models. Therefore, a consumer purchase decision model that includes successive multiple stages was considered in this study as an evolutionary game model due to its more systematic and complete nature.

The current state of research in evolutionary game theory covers a wide range of fields such as computational modelling \cite{WOS:000804643000014,WOS:000838570000046}, biology \cite{WOS:000884382400002, WOS:000954433200008} and mathematics \cite{WOS:001120954700001,WOS:000860498000011}. Some scholars have analysed different strategies and behavioural patterns based on a series of basic models from a theoretical perspective. For example, \cite{WOS:000826302800002} explored the evolutionary performance of multi-player zero determinant strategies in the Prisoner's Dilemma model. Vlastimil and Ross constructed an asymmetric Hawk-Dove game model considering time cost to discuss mixed strategies in evolutionary outcomes \cite{WOS:000841782100001}. The Prisoner's Dilemma model and the Snowdrift game model also were discussed on complex networks to validate the prediction of evolutionary outcomes by a dynamic approximation of the master equation method \cite{WOS:001120954700001}. It is worth noting that these theoretical studies have some degree of limitations. They have either explored the effects of different mechanisms based on existing classical game models, or added constraints to further refine and analyse these underlying models. In fact, one possible problem is that the classical game models do not seem to be sufficient to fully explain all the real-world subject interaction laws. Therefore, it is still of high theoretical value to try to propose some new fundamental game models to explain specific social dilemmas.

As a mathematical tool, Min-Max normalisation has been mainly applied to data processing in the field of artificial intelligence, such as machine learning, computer vision and data mining. For example, \cite{WOS:000509586800024} extended the method in the proposed adaptive activation function selection mechanism for forecasting the foreign exchange market. \cite{WOS:000673538600001} used fuzzy clustering and Min-Max normalisation to create an associative knowledge graph of video content. \cite{WOS:000886932000013} combined Min -Max normalisation into pattern recognition for invariant electromyography. In addition some relatively few theoretical researchers have tried to incorporate the tool in order to come up with new mathematical approaches. \cite{WOS:001131430800001} assigned feature weighting in Min-Max normalisation to improve the accuracy of classification methods. In conclusion, Min-Max normalisation methods are currently focused on artificial intelligence, and there is still some research gap in its application in areas such as mathematical modelling. And traditional quantitative methods directly involve factors with multi-dimensional units in the computation of evolutionary game models, which affects the computational efficiency and accuracy of the model Therefore, this study will introduce the Min-Max normalisation method in Min-Max normalisation is introduced into evolutionary game models to explore the impact of the method on the accuracy of mathematical model construction and simulation.

In summary, the innovations of our study are:

\begin{enumerate}[\quad (1)]	
	\item We proposed an expectation supply-demand game model as the underlying framework.
	\item We constructed a consumer decision-making process with feedback in the evolutionary game model.
	\item We processed the simulation values of the evolutionary game model through the Min-Max normalisation method to try to improve the accuracy of the results of the study.
	\item The study replicated the development path of NEVs and made predictions for the future.
\end{enumerate}	

{
	\linespread{1.07}
	\begin{table*}[h]
		\centering
		\vspace{-5pt}
		\begin{tabular}{|ccc|}
			\Xhline{1.2pt}
			Parameters &  Descriptions & Parameter value range  \\
			\Xhline{1.2pt}
			$x$ & Percentage of NEV manufacturers & $0 \leq x \leq 1$ \\
			$1-x$ & Percentage of TFV manufacturers & $0 \leq 1-x \leq 1$ \\
			$y$ & Percentage of NEV consumers & $0 \leq y \leq 1$ \\
			$1-y$ & Percentage of TFV consumers & $0 \leq 1-y \leq 1$ \\
			$V_1, V_2$ & Profit of NEVs and TFVs & $V_1 \geq 0, V_2 \geq 0 $ \\
			$C$ & Research and development costs of NEVs & $C \geq 0 $ \\
			$ R $ & Points bonus for NEVs in the dual credit policy & $p\geq 0$ \\
			$f_1$ & Points penalty for average fuel consumption in the dual credit policy & $f_1\geq 0$ \\
			$f_2$ & Financial penalties for causing environmental pollution & $f_1\geq 0$ \\
			$ T $ & Commuting need & $\alpha\geq 0$ \\
			$ E $ & Need for environmental awareness & $E \geq 0$ \\		
			$\alpha $ & External stimulus to demand from promotional activities and advertising & $0 \leq \alpha \leq 1$ \\
			$I_1$ & Information retrieval for vehicle offline brick-and-mortar stores & $I_1 \geq 0 $ \\
			$I_2$ & Vehicle information checking within the Internet &  $I_2 \geq 0 $ \\
			$I_3$ & Information retrieval of contacts & $I_3 \geq 0 $ \\
			$I_4$ & Retrieval of vehicle information released by the media &  $I_4 \geq 0 $ \\
			$\delta$ & Feedback factors for post-purchase behavior & $\delta \geq 0 $ \\
			$P_1, P_2$ & Prices of NEVs and TFVs&  $P_1 \geq 0, P_2 \geq 0 $ \\
			$ A $ & Vehicle purchase tax&  $ A \geq 0 $ \\
			$ p $ & Energy prices for TFVs & $p\geq 0$ \\
			$ r $ & NEV battery payout and replacement insurance proceeds & $p\geq 0$ \\
			$e_1, e_2$ & Range of NEVs and TFVs & $e_1 \geq 0, e_2 \geq 0 $ \\
			$n_1, n_2$ & Number of infrastructures for NEVs and TFVs & $n_1 \geq 0, n_2 \geq 0 $ \\
			$c_1, c_2$ & Energy supplement efficiency of NEVs and TFVs  & $c_1 \geq 0, c_2 \geq 0 $ \\ 
			\Xhline{1.2pt}	
		\end{tabular}
		\caption{\newline Descriptions of the parameters in the two-party game model.}
		\label{tab5-1}
	\end{table*}
}

\section{Methods}
\label{section3}

\subsection{Min-Max Normalisation}

Min-Max normalisation is a common method of pre-processing data, which is also known as deviation normalisation. Bylinear transformation of the original data, the method is able to scale numerical features to map into the range $[0, 1]$ or $[-1, 1]$. The transformation function is as follows:
\begin{equation}
	x_{norm} = \frac{x-x_{min}}{x_{max}-x_{min}},
\end{equation}

where, $x$ is the original data, $x_{norm}$ is the normalised data, $x_{max}$ is the maximum value of the sample data and $x_{min}$ is the minimum value. When the desired output is in the range $[-1, 1]$, the function can be: 
\begin{equation}
	X_{norm} = \frac{2*(X-X_{min})}{X_{max}-X_{min}}-1.
\end{equation}

Through the Min-Max normalisation process, the dimension expression of the original data can be transformed into a dimensionless expression, and data indicators in the same order of magnitude are more suitable for comprehensive comparison and weighting, and the calculation of the model will be simplified.

\begin{table*}[h]
	\centering 
	%\raggedright
	%\setlength{\abovecaptionskip}{0cm} % 调整间距
	%\setlength{\belowcaptionskip}{-0.9cm}
	\vspace{-5pt}
	\setlength\tabcolsep{4pt}
	\caption{ \newline Player payoff matrix for expected supply and demand games.}
	%\caption{双方博弈模型中的参数及参数说明}
	\label{tab}
	\begin{tabular}{m{3cm}m{3cm}m{3.4cm}m{3cm}}
		%\specialrule{0.8pt}{0pt}{0pt}
		%\vspace{1cm}
		\hline
		\rule{0pt}{11pt}
		& &\makecell[l]{Demand players} &  \\
		\cmidrule[0.8pt](lr){3-4}
		\rule{0pt}{11pt}
		& &\makecell[l]{cooperation (C) } & \makecell[l]{defection (D)}\\
		%\cmidrule[1.2pt](lr){2-3}
		\hline
		\rule{0pt}{17pt}
		\makecell[c]{Supply players} &\makecell[c]{ cooperation (C) } & \makecell[l]{ $\gamma + \varepsilon$, $b $} & \makecell[l]{$\gamma - \bigtriangleup \varepsilon$, $b $}\\
		%\cmidrule[0.8pt](lr){2-4}
		%\Xhline{0.8pt}
		\rule{0pt}{17pt}
		&\makecell[c]{defection (D)} & \makecell[l]{$\gamma - \bigtriangleup \varepsilon$, $b $}&
		\makecell[l]{$\gamma +  \varepsilon$, $b $}\\
		\hline
	\end{tabular}
\end{table*}

\subsection{Expected supply and demand games}

In a multi-player cooperative game, some of the players are the group $m$ that satisfies the desired demand, called the supply players, whose number is $N_m \in N$. The rest of the players are the group $\theta$ that makes the desired demand, called the demand players, whose number is $N_{\theta} \in N$. The available choice strategies for the two groups of players are cooperation (C) and defection (D), i.e., the set $\{C, D\}$. When the supply player group $m$ plays the game in tandem with the demand player group $\theta$, the payoff of supply player is:

\begin{equation}
	P_m(t)=\sum_{i\in N_m,\, j\in N_{\theta }} P_{i,j}(t),
\end{equation}      

\begin{equation}
	P_{i,j}(t)=\begin{cases}
		\gamma + \varepsilon, \,if\, S_i(t)=S_j(t)  
		\\
		\gamma - \bigtriangleup \varepsilon, \,if\, S_i(t)\neq S_j(t)
	\end{cases}
\end{equation}
where $P_m(t)$ is the total payoff of the supply players, $P_{i,j}(t)$ is the payoff of the game between individual player $i$ in the supply player group and individual player $j$ in the demand player group at moment $t$, and $S_i(t)$ and $S_j(t)$ are the strategies of $i$ and $j$, respectively. $\gamma$ is the basic payoff of supply players, $\varepsilon$ is the expected matching payoff when the two players' strategies are consistent, and when their strategy choices are inconsistent then the demand players with dissatisfaction generates an expected mismatch discount $\bigtriangleup \varepsilon = \delta \varepsilon$, $\delta$ is a discount factor and $\delta \in [0,\,1)$. Similarly, the total payoff of the demand players is:

\begin{equation}
	P_{\theta}(t)=\sum_{j\in N_{\theta}, \, i\in N_m } P_{j,i}(t).
\end{equation}
The model assumes that the psychological gains of the demand players are not taken into account at this point. The supply players in the usual case cannot impose an expected discount on the demand players, so the demand player's payoff fails to change. The payoff matrix for the expected supply and demand game is shown in Table \ref{tab}.

\subsection{Basic assumptions of evolutionary game model framework}

In this two-sided evolutionary game model, the vehicle manufacturer's alternative strategies are to produce NEVs ($x$) and to produce TFVs ($1-x$). Consumers' strategies are to buy NEVs ($y$) and to buy TFVs ($1-y$). For vehicle manufacturers, the study supposed that the NEV point bonus they receive in the dual-credit policy is denoted as $R$ and the average fuel consumption point penalty is denoted as $f_1$. Manufacturers can benefit $V_1$ from producing NEVs, but they must also pay $C_1$ for research and development, but they also need to bear the research and development costs $C_1$. Similarly, the manufacturers who choose to produce TFVs can make a profit $V_2$ and receive a pollution penalty $f_2$.

%[scale=0.4] [width=10cm,height=5cm]

\begin{figure*}[h]
	\subfigure{	
		\begin{minipage}[]{1\linewidth}
			\centering
			\includegraphics[width=11cm,height=5.5cm]{"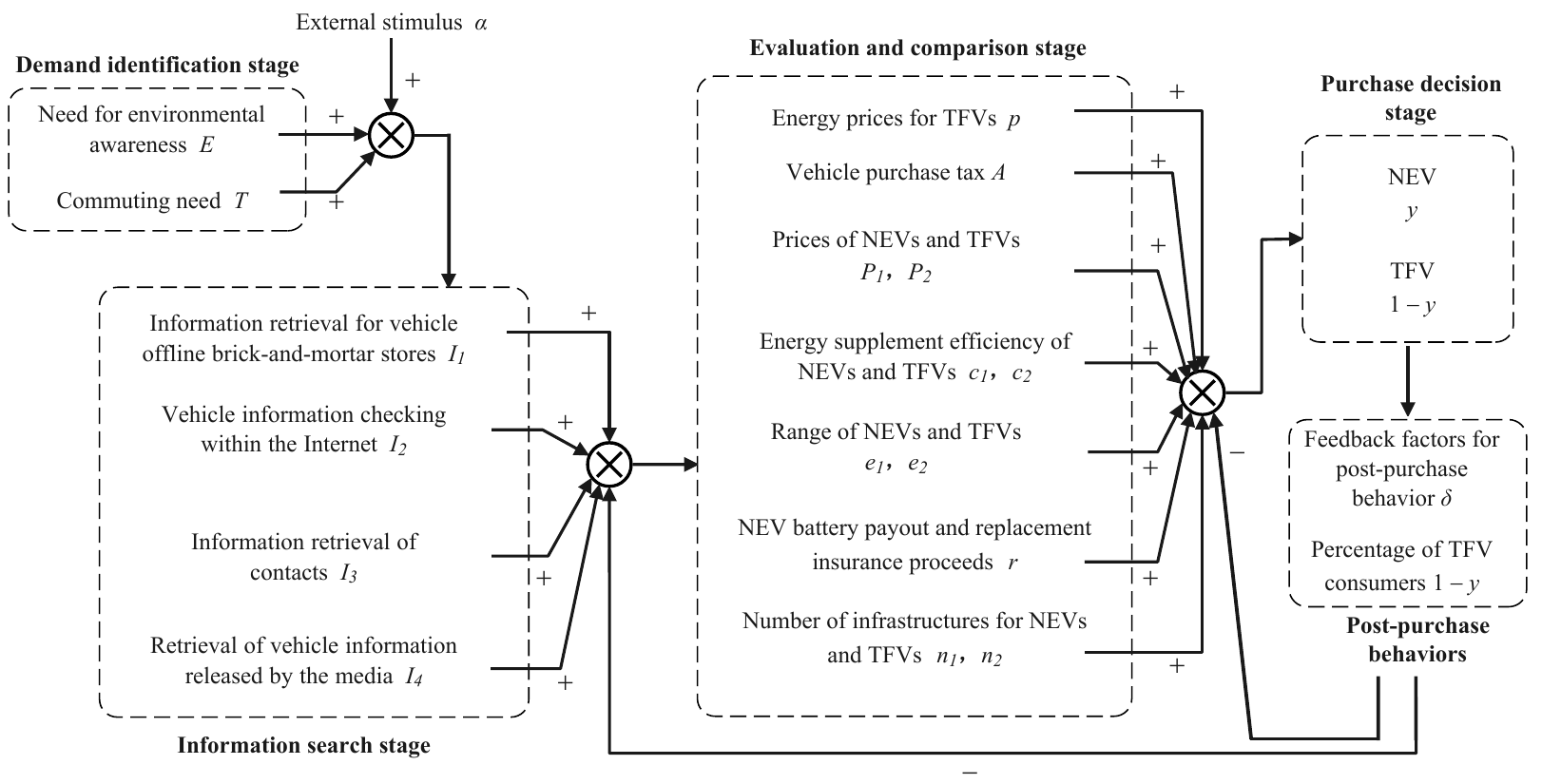"}
		\end{minipage}
	}
	\vspace{-5pt}
	\caption{\centering Consumer decision-making process with feedback for NEVs.}
	\label{FIG5-1}
\end{figure*}

For the consumer level, we introduced the feedback factor of expectation matching from the proposed supply and demand game in the post-purchase behavior, and proposed the consumer decision-making process \cite{0887302X8700600104,mar20322} of NEVs with feedback. As shown in Fig.\ref{FIG5-1}, in the demand identification stage consumers can receive external stimuli $\alpha$ such as promotions or advertisements to elicit their needs or problems, which include basic commuting needs $T$ with environmentally conscious needs $E$. In the information search stage, consumers will try to obtain information about relevant feasible solutions to realize their needs. We assumeed that their main external sources of information collection are: NEVs through offline brick-and-mortar stores $I_1$, NEV information on the Internet $I_2$, contacts and inquiries $I_3$, and NEV information released by the media $I_4$. After this, consumers consider various factors to evaluate and compare feasible solutions, such as fuel price $p$, vehicle purchase tax $A$, price of NEVs $P_1$ and TFVs $P_2$, efficiency of energy supplementation $c_1$, $ c_2$, range $e_1$, $ e_2$, and amount of infrastructure $n_1$, $n_2$. Meanwhile, we also introduced payout and replacement insurance$r$ for NEV batteries in this stage to explore feasible incentive mechanisms. Through these stages, consumers will make a purchase decision to buy NEVs or TFVs. When consumers experienced the vehicle, they will generate a series of post-purchase behaviors such as re-purchase, complaint, and sharing. A feedback factor $\sigma$ is set to simulate the impact of their post-purchase behavior on the reputation of the vehicle manufacturer, which may ultimately lead to a loss of benefit in the information-seeking stage. The percentage of consumers of TFVs will also affect the evaluation and comparison stage as a feedback of post-purchase behavior. Related symbols and definitions are listed in Table \ref{tab5-1}.

\subsection{System dynamics model}
\begin{figure*}[h!]
	\subfigure{	
		\begin{minipage}[]{1\linewidth}
			\centering
			\includegraphics[scale=0.43]{"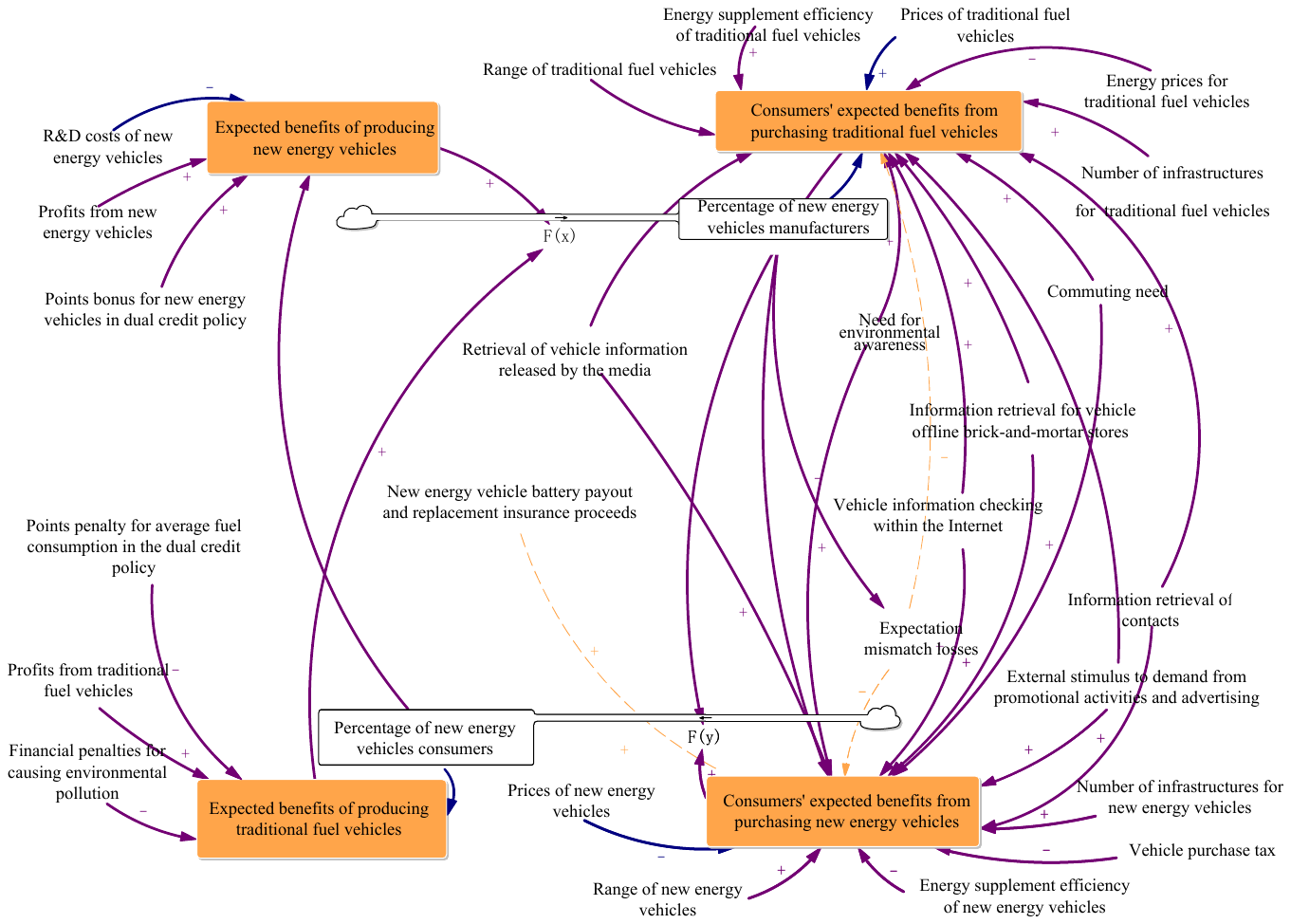"}
		\end{minipage}
	}
	\vspace{-5pt}
	\caption{\centering System dynamics model without feedback.}
	\label{FIG5-2}
\end{figure*}

\begin{figure*}[h!]
	\subfigure{	
		\begin{minipage}[]{1\linewidth}
			\centering
			\includegraphics[scale=0.43]{"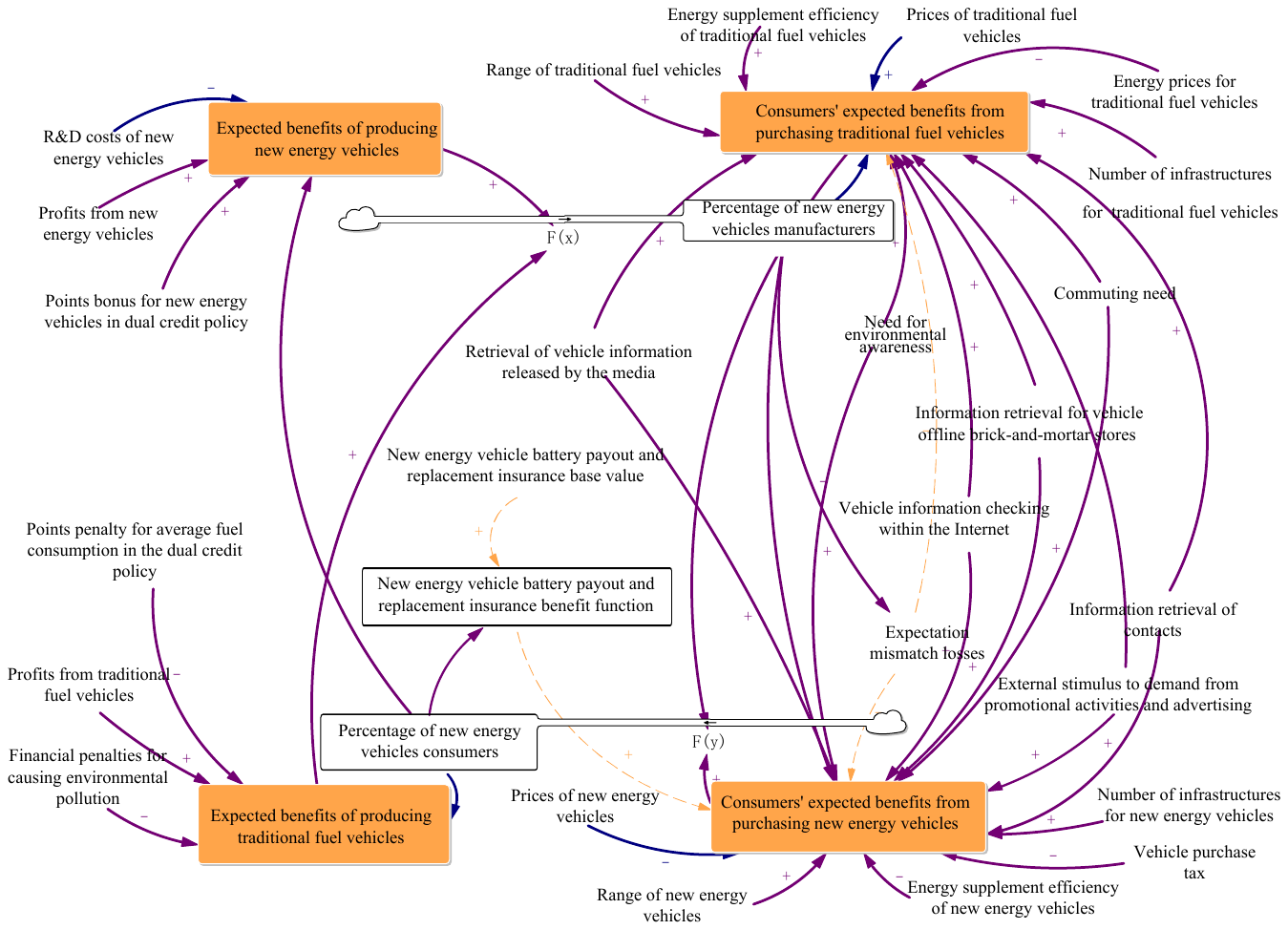"}
		\end{minipage}
	}
	\vspace{-5pt}
	\caption{\centering System dynamics model with feedback.}
	\label{FIG5-3}
\end{figure*}

As shown in Fig.\ref{FIG5-2} and Fig.\ref{FIG5-3}, based on the assumptions, the study modeled the system dynamics between the vehicle manufacturer and the consumer. The system dynamics model without feedback includes 2 flow rates, 2 rate variables, 24 macro variables and 5 medium variables. The system dynamics model with feedback has 1 more intermediate variable and 1 more macro variable. 

%\vspace{-5pt}
\section{Model analysis}
\label{section4}

\subsection{Replication dynamic analysis of	the model}

{
	\linespread{1.07}
	\begin{table*}[h!]
		\centering 
		%\raggedright
		%\setlength{\abovecaptionskip}{0cm} % 调整间距
		%\setlength{\belowcaptionskip}{-0.9cm}
		\vspace{-5pt}
		\setlength\tabcolsep{4pt}
		\caption{ \newline Payoff matrix among agents.}
		\begin{tabular}{m{3cm}m{3cm}m{3.4cm}m{3cm}}
			%\specialrule{0.8pt}{0pt}{0pt}
			%\vspace{1cm}
			\hline
			\rule{0pt}{11pt}
			& &\makecell[l]{Consumers} &  \\
			\cmidrule[0.8pt](lr){3-4}
			\rule{0pt}{11pt}
			& &\makecell[l]{NEVs\,$(y)$} & \makecell[l]{TFVs\,$(1-y)$}\\
			%\cmidrule[1.2pt](lr){2-3}
			\hline
			\rule{0pt}{17pt}
			\makecell[c]{Vehicle manufacturer} &\makecell[c]{NEVs\,$(x)$} & \makecell[l]{ $\Upsilon_1$, $ \Upsilon_3 $} & \makecell[l]{$\Upsilon_1$, $\Upsilon_4$}\\
			%\cmidrule[0.8pt](lr){2-4}
			%\Xhline{0.8pt}
			\rule{0pt}{17pt}
			&\makecell[c]{TFVs\,$(1-x)$} & \makecell[l]{$\Upsilon_2$,
				$\Upsilon_3$}&
			\makecell[l]{$\Upsilon_2$,
				$\Upsilon_4$}\\
			\hline
		\end{tabular}
		\label{tab5-2}
	\end{table*}
}

As shown in Table \ref{tab5-2} are the payoff matrices for the evolutionary game model between vehicle manufacturers and consumers without the introduction of feedback. The payoffs are as follows:
\begin{eqnarray}
	\Upsilon_1 &=& R+V_1-C, \nonumber \\
	\Upsilon_2 &=& V_2-f_1-f_2, \nonumber \\
	\Upsilon_3 &=& \alpha_1 (T+E)+\varepsilon(I_1+I_2+I_3+I_4)+P_1+e_1+n_1-c_1+r-A, \nonumber \\
	\Upsilon_4 &=& \alpha_1 (T+E)+ \varepsilon(I_1+I_2+I_3+I_4)+P_2+e_2+n_2-c_2-p. \nonumber 
\end{eqnarray}

Accordingly, the expected utility of vehicle manufacturers to produce NEVs is $U_{X1}$:
\begin{eqnarray}
	U_{X1} &=& R+V_1-C. 
\end{eqnarray}

The expected utility $U_{X2}$ of producing TFVs is:
\begin{eqnarray}
	U_{X2} &=& V_2-f_1-f_2.  
\end{eqnarray}

The average expected utility $\overline U_X$ is:
\begin{eqnarray}
	\overline U_{X} &=& xU_{X1}+(1-x)U_{X2}.  
\end{eqnarray}

Then the replicated dynamic equation $F(x)$ can be obtained \cite{Evolutionary}:  
\begin{eqnarray}
	\label{5eq:fx}
	F(x) =\frac{dx}{dt}=x(U_{X1} - \overline U_X) = x (1-x) [(R+V_1-C)-(V_2-f_1-f_2)].
\end{eqnarray}

The derivative of $F(x)$ yields $F'(x)$ is: 
\begin{eqnarray}
	\label{5eq:f'x}
	{F}'(x) = \frac{dF(x)}{dx}= (1-2x)[(R+V_1-C)-(V_2-f_1-f_2)].
\end{eqnarray}

\begin{figure*}[h]
	\subfigure[$R+V_1-C \textgreater V_2-f_1-f_2$]{\label{FIG5-4a}	
		\begin{minipage}[]{.45\linewidth}
			\centering
			\includegraphics[scale=0.15]{"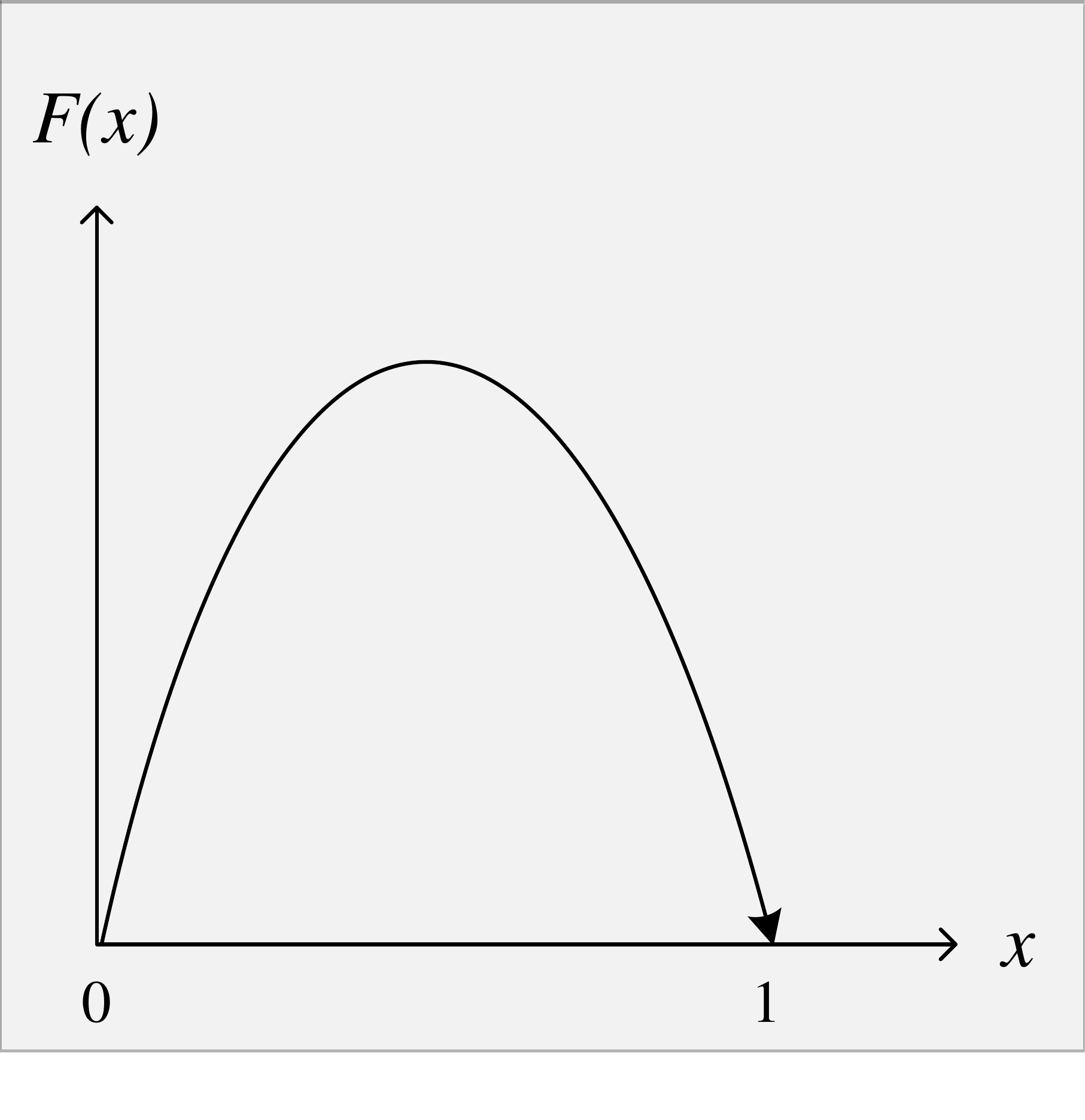"}
	\end{minipage}	}
	\subfigure[$R+V_1-C  \textless V_2-f_1-f_2$]{\label{FIG5-4b}
		\begin{minipage}[]{.45\linewidth}
			\centering
			\includegraphics[scale=0.15]{"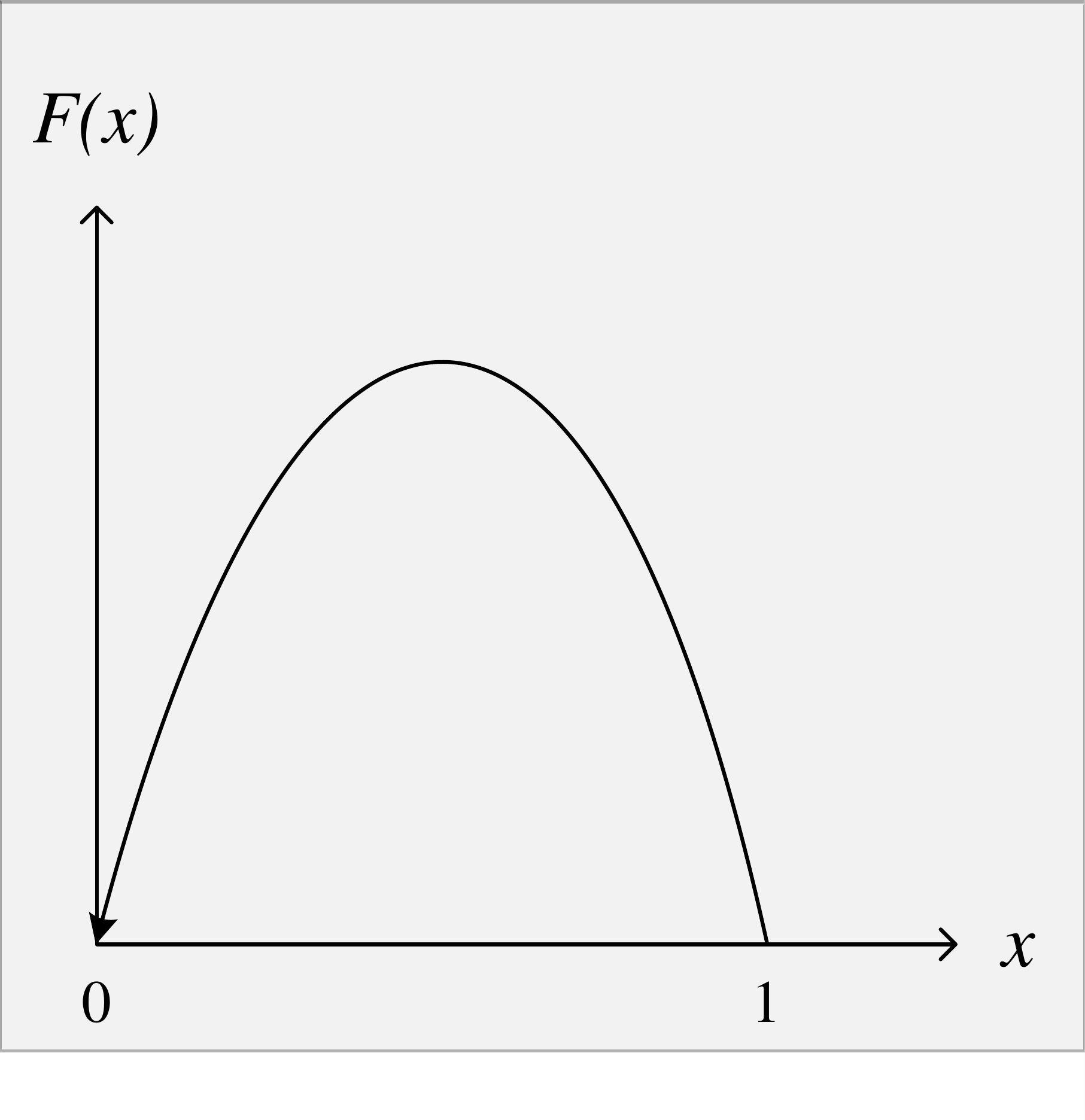"}
		\end{minipage}
	}
	%\vspace{-5pt}
	\caption{\centering The phase portrait of the vehicle manufacturer.}
\end{figure*}

Solving for $F(x) = 0$ yields  $x=0$ and $x=1$. The two cases need to be discussed separately according to the formula \eqref{5eq:f'x}:

(1) If $R+V_1-C \textgreater V_2-f_1-f_2$, we can obtain: ${F}'(x)|_{x=1} \textless 0 $ and ${F}'(x)|_{x=0} \textgreater 0 $. Then, as shown in Fig.\ref{FIG5-4a}, $x=1$ is only ESS.

(2) If $R+V_1-C  \textless V_2-f_1-f_2$, we can obtain: ${F}'(x)|_{x=0} \textless 0 $ and ${F}'(x)|_{x=1} \textgreater 0 $. Then, as shown in Fig.\ref{FIG5-4b}, $x=0$ is only ESS.

Similarly, the expected utility of consumers buying NEVs $U_{Y1}$ is:
\begin{eqnarray}
	U_{Y1} &=& \alpha_1 (T+E)+\varepsilon(I_1+I_2+I_3+I_4)+P_1+e_1+n_1-c_1+r-A. 
\end{eqnarray}

The expected utility of consumers buying TFVs $U_{Y2}$is:
\begin{eqnarray}
	U_{Y2} &=& \alpha_1 (T+E)+ \varepsilon(I_1+I_2+I_3+I_4)+P_2+e_2+n_2-c_2-p.  
\end{eqnarray}

The average expected utility $\overline U_Y$ is:
\begin{eqnarray}
	\overline U_{Y} &=& yU_{Y1}+(1-y)U_{Y2}.  
\end{eqnarray}

The consumer's replication dynamic equation $F(y)$ is:
\begin{eqnarray}
	\label{5eq:fy}
	F(y) &=&\frac{dy}{dt}=x(U_{Y1} - \overline U_Y) \nonumber \\
	&=& y (1-y) [(P_1+e_1+n_1-c_1+r-A)-(P_2+e_2+n_2-c_2-P)].
\end{eqnarray}

The derivative of $F(y)$  yields $F'(y)$ is: 
\begin{eqnarray}
	\label{5eq:f'y}
	{F}'(y) &=& \frac{dF(y)}{dy} \nonumber \\
	&=& (1-2y)[(P_1+e_1+n_1-c_1+r-A)-(P_2+e_2+n_2-c_2-P)].
\end{eqnarray}

\begin{figure*}[h!]
	\subfigure[$P_1+e_1+n_1-c_1+r-A \textgreater P_2+e_2+n_2-c_2-P$]{\label{FIG5-5a}	
		\begin{minipage}[]{.45\linewidth}
			\centering
			\includegraphics[scale=0.15]{"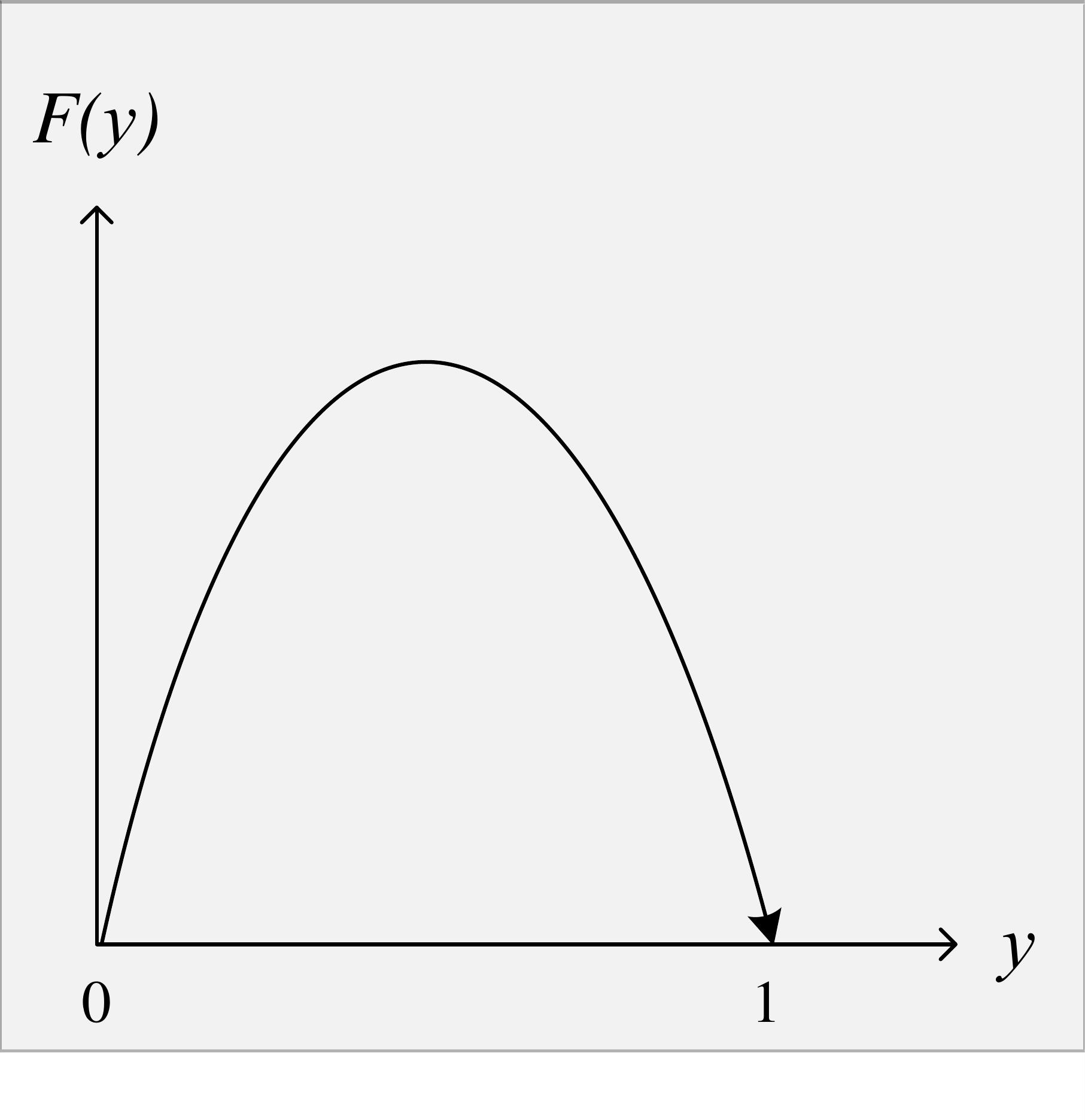"}
	\end{minipage}	}
	\subfigure[$P_1+e_1+n_1-c_1+r-A \textless P_2+e_2+n_2-c_2-P$]{\label{FIG5-5b}
		\begin{minipage}[]{.45\linewidth}
			\centering
			\includegraphics[scale=0.15]{"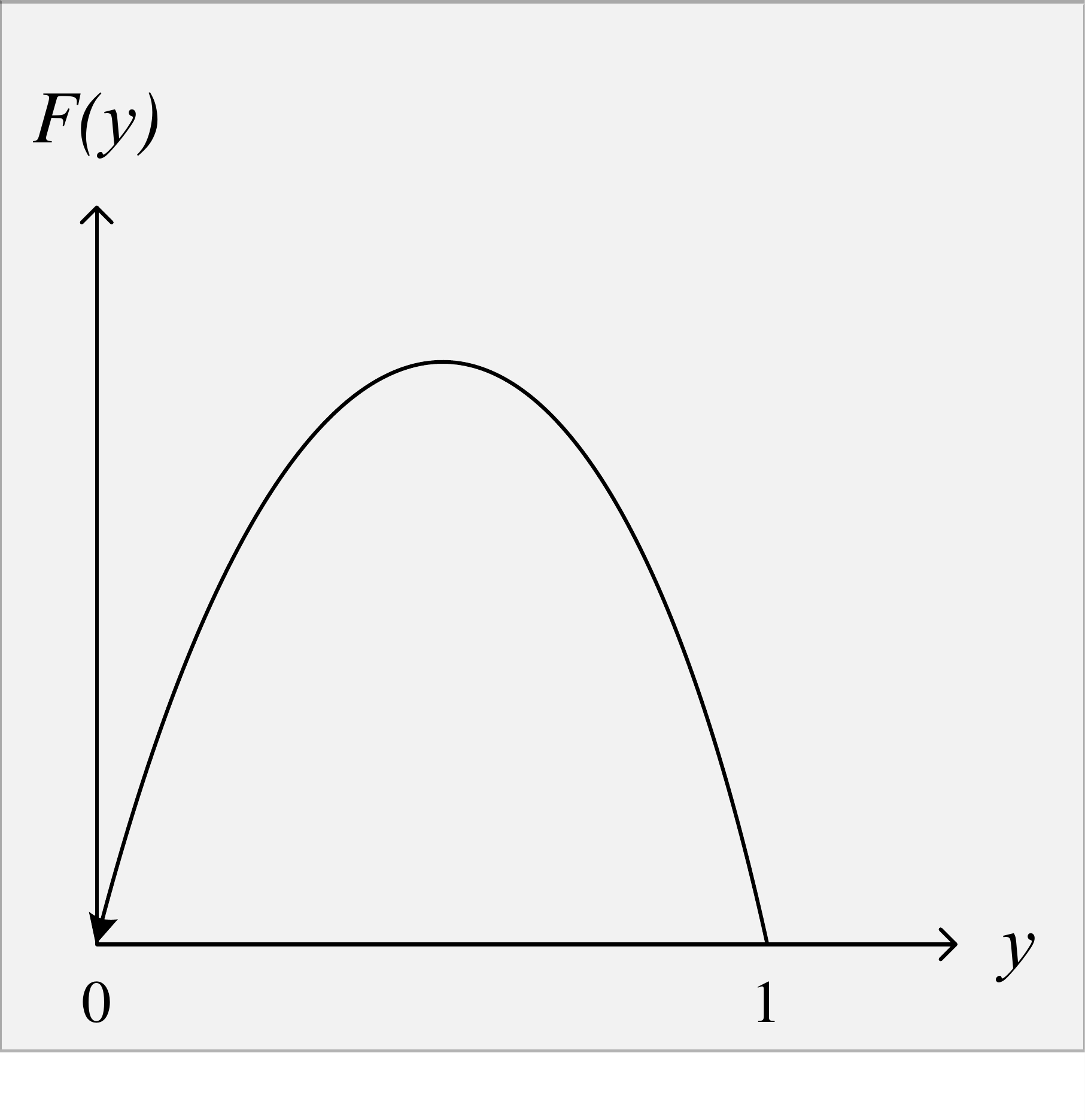"}
		\end{minipage}
	}
	%\vspace{-5pt}
	\caption{\centering The phase portrait of the customer.}
\end{figure*}

Solving for $F(y)=0 $ yields $y=0$ and $y=1$. Two separate cases need to be discussed based on the equation \eqref{5eq:f'y}:

(1) If  $ P_1+e_1+n_1-c_1+r-A \textgreater P_2+e_2+n_2-c_2-P$, we can obtain:  ${F}'(y)|_{y=1} \textless 0 $ and ${F}'(y)|_{y=0} \textgreater 0 $. Then, as shown in Fig.\ref{FIG5-5a}, $y=1$ is only ESS.

(2) If $P_1+e_1+n_1-c_1+r-A \textless P_2+e_2+n_2-c_2-P$, we can obtain: ${F}'(y)|_{y=0} \textless 0 $ and ${F}'(y)|_{y=1} \textgreater 0 $. Then, as shown in Fig.\ref{FIG5-5b}, $y=0$ is only ESS.

\subsection{Stability analysis of systems without feedback}

The Jacobi matrix $J_{I}$ is:
\begin{eqnarray}
	\label{5JI}
	J_{I}=
	\begin{bmatrix}
		\frac{\partial F(x)}{\partial x} &  \frac{\partial F(x)}{\partial y} \vspace{5pt}\\
		\frac{\partial F(y)}{\partial x} &  \frac{\partial F(y)}{\partial y}
	\end{bmatrix}
	=	
	\begin{bmatrix}
		\Omega _1  & 0 \vspace{5pt}  \\
		0 & \Omega _2
	\end{bmatrix}.
\end{eqnarray}

where,
\begin{align*}
	\Omega _1=&(1-2x)[(R+V_1-C)-(V_2-f_1-f_2)],\\
	\Omega _2=&(1-2y)[(P_1+e_1+n_1-c_1+r-A)-(P_2+e_2+n_2-c_2-P)].
\end{align*}

The equilibrium points of the system without feedback are: $(\,0\, ,\,0\,)$, $(\,0\, ,\,1\,)$, $(\,1\, ,\,0\,)$, $(\,1\, ,\,1\,)$.

\textbf{Proposition 1}.\quad Each of the points $(\,1\,,\,0\,)$, $(\,1\,,\,1\,)$, $(\,0\,,\,0\,)$, and $(\,0\,,\,1\,)$ can be saddle or unstable in the model without feedback, which is only characterised by the emergence of stable points.

\textbf{Proof 1}.\quad The ESS points of the system require $det(J_I)>0$ and $tr(J_I)<0$.

(1) When $R+V_1-C \textgreater V_2-f_1-f_2$ and   $P_1+e_1+n_1-c_1+r-A \textgreater P_2+e_2+n_2-c_2-P$, each equilibrium is analysed as shown in Table \ref{tab5-4}.

(2) When $R+V_1-C \textgreater V_2-f_1-f_2$ and  $P_1+e_1+n_1-c_1+r-A \textless P_2+e_2+n_2-c_2-P$, each equilibrium is analysed as shown in Table \ref{tab5-5}.

(3) When $R+V_1-C \textless V_2-f_1-f_2$ and  $P_1+e_1+n_1-c_1+r-A \textgreater P_2+e_2+n_2-c_2-P$, each equilibrium is analysed as shown in Table \ref{tab5-6}.

(4) When $R+V_1-C \textless V_2-f_1-f_2$ and  $P_1+e_1+n_1-c_1+r-A \textless P_2+e_2+n_2-c_2-P$, each equilibrium is analysed as shown in Table \ref{tab5-7}.

{
	\linespread{1.07}
	\begin{table*}[h]
		\centering 
		%\raggedright
		%\setlength{\abovecaptionskip}{0cm} % 调整间距
		%\setlength{\belowcaptionskip}{-0.9cm}
		\vspace{-5pt}
		\setlength\tabcolsep{4pt}
		\caption{\newline Stability analysis of the two-player model in case one.}
		\label{tab5-4}
		\begin{tabular}{cccc}
			\hline
			\rule{0pt}{12pt}
			Equilibrium point $(\,x\,,\,y\,)$ & \qquad Determinant$(J_{I})$ & \qquad Traces$(J_{I})$ & \qquad Results\\
			\hline
			\makecell[c]{ $(\,0\,,\,0\,)$} & 
			\makecell[c]{ \qquad +}& 
			\makecell[c]{ \qquad +}& 
			\makecell[c]{ \qquad Unstable point}\\
			
			\makecell[c]{ $(\,0\,,\,1\,)$} & 
			\makecell[c]{ \qquad -}& 
			\makecell[c]{ \qquad Uncertain}& 
			\makecell[c]{ \qquad Saddle point}\\
			
			\makecell[c]{ $(\,1\,,\,0\,)$} & 
			\makecell[c]{ \qquad -}& 
			\makecell[c]{ \qquad Uncertain}& 
			\makecell[c]{ \qquad Saddle point}\\
			
			\makecell[c]{ $(\,1\,,\,1\,)$} & 
			\makecell[c]{ \qquad +}& 
			\makecell[c]{ \qquad -}& 
			\makecell[c]{ \qquad ESS}\\
			\hline	
		\end{tabular}
	\end{table*}
}

{
	\linespread{1.07}
	\begin{table*}[h]
		\centering 
		%\raggedright
		%\setlength{\abovecaptionskip}{0cm} % 调整间距
		%\setlength{\belowcaptionskip}{-0.9cm}
		\vspace{-5pt}
		\setlength\tabcolsep{4pt}
		\caption{\newline Stability analysis of  the two-player model in case two.}
		\label{tab5-5}
		\begin{tabular}{cccc}
			\hline
			\rule{0pt}{12pt}
			Equilibrium point $(\,x\,,\,y\,)$ & \qquad Determinant$(J_{I})$ & \qquad Traces$(J_{I})$ & \qquad Results\\
			\hline
			\makecell[c]{ $(\,0\,,\,0\,)$} & 
			\makecell[c]{ \qquad -}& 
			\makecell[c]{ \qquad Uncertain}& 
			\makecell[c]{ \qquad Saddle point}\\
			
			\makecell[c]{ $(\,0\,,\,1\,)$} & 
			\makecell[c]{ \qquad +}& 
			\makecell[c]{ \qquad +}& 
			\makecell[c]{ \qquad Unstable point}\\
			
			\makecell[c]{ $(\,1\,,\,0\,)$} & 
			\makecell[c]{ \qquad +}& 
			\makecell[c]{ \qquad -}& 
			\makecell[c]{ \qquad ESS}\\
			
			\makecell[c]{ $(\,1\,,\,1\,)$} & 
			\makecell[c]{ \qquad -}& 
			\makecell[c]{ \qquad Uncertain}& 
			\makecell[c]{ \qquad Saddle point}\\
			\hline	
		\end{tabular}
	\end{table*}
}

{
	\linespread{1.07}
	\begin{table*}[h]
		\centering 
		%\raggedright
		%\setlength{\abovecaptionskip}{0cm} % 调整间距
		%\setlength{\belowcaptionskip}{-0.9cm}
		\vspace{-5pt}
		\setlength\tabcolsep{4pt}
		\caption{\newline Stability analysis of the two-player model in case three.}
		\label{tab5-6}
		\begin{tabular}{cccc}
			\hline
			\rule{0pt}{12pt}
			Equilibrium point $(\,x\,,\,y\,)$ & \qquad Determinant$(J_{I})$ & \qquad Traces$(J_{I})$ & \qquad Results\\
			\hline
			\makecell[c]{ $(\,0\,,\,0\,)$} & 
			\makecell[c]{ \qquad -}& 
			\makecell[c]{ \qquad Uncertain}& 
			\makecell[c]{ \qquad Saddle point}\\
			
			\makecell[c]{ $(\,0\,,\,1\,)$} & 
			\makecell[c]{ \qquad +}& 
			\makecell[c]{ \qquad -}& 
			\makecell[c]{ \qquad ESS}\\
			
			\makecell[c]{ $(\,1\,,\,0\,)$} & 
			\makecell[c]{ \qquad +}& 
			\makecell[c]{ \qquad +}& 
			\makecell[c]{ \qquad Unstable point}\\
			
			\makecell[c]{ $(\,1\,,\,1\,)$} & 
			\makecell[c]{ \qquad -}& 
			\makecell[c]{ \qquad Uncertain}& 
			\makecell[c]{ \qquad Saddle point}\\
			\hline	
		\end{tabular}
	\end{table*}
	
}

{
	\linespread{1.07}
	\begin{table*}[h]
		\centering 
		%\raggedright
		%\setlength{\abovecaptionskip}{0cm} % 调整间距
		%\setlength{\belowcaptionskip}{-0.9cm}
		\vspace{-5pt}
		\setlength\tabcolsep{4pt}
		\caption{\newline Stability analysis of  the two-player model in case four.}
		\label{tab5-7}
		\begin{tabular}{cccc}
			\hline
			\rule{0pt}{12pt}
			Equilibrium point $(\,x\,,\,y\,)$ & \qquad Determinant$(J_{I})$ & \qquad Traces$(J_{I})$ & \qquad Results\\
			\hline
			\makecell[c]{ $(\,0\,,\,0\,)$} & 
			\makecell[c]{ \qquad +}& 
			\makecell[c]{ \qquad -}& 
			\makecell[c]{ \qquad ESS}\\
			
			\makecell[c]{ $(\,0\,,\,1\,)$} & 
			\makecell[c]{ \qquad -}& 
			\makecell[c]{ \qquad Uncertain}& 
			\makecell[c]{ \qquad Saddle point}\\
			
			\makecell[c]{ $(\,1\,,\,0\,)$} & 
			\makecell[c]{ \qquad -}& 
			\makecell[c]{ \qquad Uncertain}& 
			\makecell[c]{ \qquad Saddle point}\\
			
			\makecell[c]{ $(\,1\,,\,1\,)$} & 
			\makecell[c]{ \qquad +}& 
			\makecell[c]{ \qquad +}& 
			\makecell[c]{ \qquad Unstable point}\\
			\hline	
		\end{tabular}
	\end{table*}

}

\subsection{Stability analysis of systems with feedback}

With the introduction of post-purchase behavioural feedback, the game model generates an expectation mismatch loss. The payout of NEV batteries and replacement insurance proceeds will vary according to the number of TFVs, with $r$ replaced by $R(y)=r(1-y)$. At this time the Jacobi matrix of the replicated dynamic system with feedback is:

\begin{eqnarray}
	\label{5JII}
	J_{I}=
	\begin{bmatrix}
		\frac{\partial F(x)}{\partial x} &  \frac{\partial F(x)}{\partial y} \vspace{5pt}\\
		\frac{\partial F(y)}{\partial x} &  \frac{\partial F(y)}{\partial y}
	\end{bmatrix}
	=	
	\begin{bmatrix}
		\Omega_3  & 0 \vspace{5pt}  \\
		\Omega_4 & \Omega_5
	\end{bmatrix}.
\end{eqnarray}

where,
\begin{align*}
	\Omega_3=&(1-2x)[(R+V_1-C)-(V_2-f_1-f_2)],\\
	\Omega_4=&2(1-2y)(\varepsilon - \bigtriangleup \varepsilon )(I_1+I_2+I_3+I_4),\\
	\Omega_5=&(1-2y)\{2x(\varepsilon - \bigtriangleup \varepsilon)(I_1+I_2+I_3+I_4)[P_1+e_1+n_1-c_1+R(y)-A] \\
	-&(P_2+e_2+n_2-c_2-P)+(\bigtriangleup \varepsilon-\varepsilon)(I_1+I_2+I_3+I_4)\}+y(1-y)R'(y).
\end{align*}

The equilibrium point of the system with feedback is $(\,0\,,\,0\,)$, $(\,0\,,\,1\,)$, $(\,1\,,\,0\,)$, $(\,x^*\,,\,y^*\,)$, where $x^* = \{ [P_1+e_1+n_1-c_1+R(y)-A]-(P_2+e_2+n_2-c_2-P)+ (\bigtriangleup \varepsilon-\varepsilon)(I_1+I_2+I_3+I_4) \} /[2(\varepsilon - \bigtriangleup \varepsilon)(I_1+I_2+I_3+I_4)]$, $y^*= [(P_1+e_1+n_1-c_1-A)-(P_2+e_2+n_2-c_2-P)+2x(\varepsilon - \bigtriangleup \varepsilon)(I_1+I_2+I_3+I_4)+(\bigtriangleup \varepsilon-\varepsilon)(I_1+I_2+I_3+I_4)+r]/r $.

\textbf{Proposition 2}.\quad Table \ref{tab5-9} shows the analysis for the stability of equilibrium points. They satisfy $0\leq \{ [P_1+e_1+n_1-c_1+R(y)-A]-(P_2+e_2+n_2-c_2-P)+ (\bigtriangleup \varepsilon-\varepsilon)(I_1+I_2+I_3+I_4) \} /[2(\varepsilon - \bigtriangleup \varepsilon)(I_1+I_2+I_3+I_4)] \leq 1 $ and $0\leq  [(P_1+e_1+n_1-c_1-A)-(P_2+e_2+n_2-c_2-P)+2x(\varepsilon - \bigtriangleup \varepsilon)(I_1+I_2+I_3+I_4)+(\bigtriangleup \varepsilon-\varepsilon)(I_1+I_2+I_3+I_4)+r]/r \leq 1$. 

(1) The saddle points of the system with feedback are $(\,0\,,0\,)$, $(\,1\,,\,0\,)$ and $(\,1\,,\,1\,)$. The point $(\,0\,,\,1\,)$ is the unstable point.

(2) The point $(\,x^o\,,\,y^o\,)$ in the system with feedback is ESS.

\textbf{Proof 2}.\quad When the $(x, y)$ is $(\,0\,,0\,)$, $(\,1\,,\,0\,)$ and $(\,1\,,\,1\,)$, the sign of $det(J_{II})$ is negative, so they are saddle points. When the $(x, y)$ is $(\,0\,,\,1\,)$, $det(J_{II})>0$ and $tr(J_{II})>0$, so this point is unstable. 

At point $(\,x^*\,,\,y^*\,)$, the Jacobi matrix $J_{II}'$ is:
\begin{eqnarray}
	\label{5JII'}
	J_{II}'=
	\begin{bmatrix}
		(1-2x^*)[(R+V_1-C)-(V_2-f_1-f_2)] & 0 \vspace{8pt}\\
		2(1-2y^*)(\varepsilon - \bigtriangleup \varepsilon )(I_1+I_2+I_3+I_4)&y^*(1-y^*)R'(y) 
	\end{bmatrix}.
\end{eqnarray}

The characteristic roots of equation $|\lambda E-J_{II}'|=0 $ is: 
$$
\lambda_{1,2}=\dfrac{(1-2x^*)[(R+V_1-C)-(V_2-f_1-f_2)]+y^*(1-y^*)R'(y)}{2} \pm \sqrt{\bigtriangleup}.
$$

The Jacobi matrix $J_{II}'$ contains the real eigenroots $\lambda_1$, $\lambda_2$ . Therefore, $(\,x^*\,,\,y^*\,)$ is the ESS of this model. The percentage of vehicle manufacturers choosing to produce NEVs is $x^*$, and the percentage of consumers buying NEVs is $y^*$.

{
	\linespread{1.07}
	
	\begin{table*}[h!]
		\centering 
		%\raggedright
		%\setlength{\abovecaptionskip}{0cm} % 调整间距
		%\setlength{\belowcaptionskip}{-0.9cm}
		\vspace{-5pt}
		\setlength\tabcolsep{4pt}
		\caption{\newline Stability analysis of replicated dynamical system with feedback. }
		\label{tab5-9}
		\begin{tabular}{cccc}
			\hline
			\rule{0pt}{12pt}
			Equilibrium point $(\,x\,,\,y\,)$ & \qquad Determinant$(J_{II})$ & \qquad Traces$(J_{II})$ & \qquad Results\\
			\Xhline{1.2pt}
			\makecell[c]{ $(\,0\,,\,0\,)$} & 
			\makecell[c]{ \qquad -}& 
			\makecell[c]{ \qquad +}& 
			\makecell[c]{ \qquad Saddle point}\\
			
			\makecell[c]{ $(\,0\,,\,1\,)$} & 
			\makecell[c]{ \qquad +}& 
			\makecell[c]{ \qquad +}& 
			\makecell[c]{ \qquad Unstable point}\\
			
			\makecell[c]{ $(\,1\,,\,0\,)$} & 
			\makecell[c]{ \qquad -}& 
			\makecell[c]{ \qquad -}& 
			\makecell[c]{ \qquad Saddle point}\\
			
			\makecell[c]{ $(\,1\,,\,1\,)$} & 
			\makecell[c]{ \qquad -}& 
			\makecell[c]{ \qquad -}& 
			\makecell[c]{ \qquad Saddle point}\\
			
			\makecell[c]{$(\,x^o\,,\,y^o\,)$} & 
			\makecell[c]{ \qquad +}& 
			\makecell[c]{ \qquad -}& 
			\makecell[c]{ \qquad ESS}\\
			\hline	
		\end{tabular}
	\end{table*}
}

\section{Simulation results}
\label{section5}

\subsection{Key parameter settings}

By the Ministry of Industry and Information Technology \cite{2021qiche} released the "December 2021 automobile industry economic operation" can be known, in 2021 the NEV production and sales were completed 3,545,000 and 3,521,000, each accounting for 13.5\% and 13.4\% of the industry, the initial value of $x $ and $y $ according to which were set to 0.135 and 0.134, respectively. Let the selling price of NEVs is $P_1=290,000$, the selling price of TFVs is $P_2=150,000$, and the vehicle purchase tax  is $A=13,274$. At the performance level, the range of NEVs and TFVs is set to $e_1=500$km and $e_2$=700km respectively, and the efficiency of energy replenishment is set to $c_1=120$min and $c_2=15$min respectively. According to \cite{2021qiche}, NEVs can generate $V_1=60.6$ billion of profit for the automakers, while TFVs can generate $V_1=60.6$ billion of profit for the automakers. In addition, TFVs can generate a profit of $V_2=3882$ billion. Traditional evolutionary game data do not pre-process quantitative data, and different parameters usually have different units of measure and ranges, and the quantitative analysis of evolutionary game models may suffer from reduced accuracy due to these inconsistencies. For this reason, this chapter attempts to reduce the effect of the range of measures by normalising the data through Min-Max, so that different factors have the same weights to improve the accuracy of the model. The values of these parameters after Min-Max normalisation are $P_1=1$, $A=0.4705$, $P_2=-1$, $e_1=-1$, $e_2=1$, $c_1=1$, $c_2=-1$, $V_1=-1$, $V_2=1$.

\subsection{The effect of feedback regulation on evolution}

\begin{figure*}[h!]
	\subfigure[Evolutionary trajectories]{\label{FIG5-6a}	
		\begin{minipage}[]{0.46\linewidth}
			\centering
			\includegraphics[scale=0.45]{"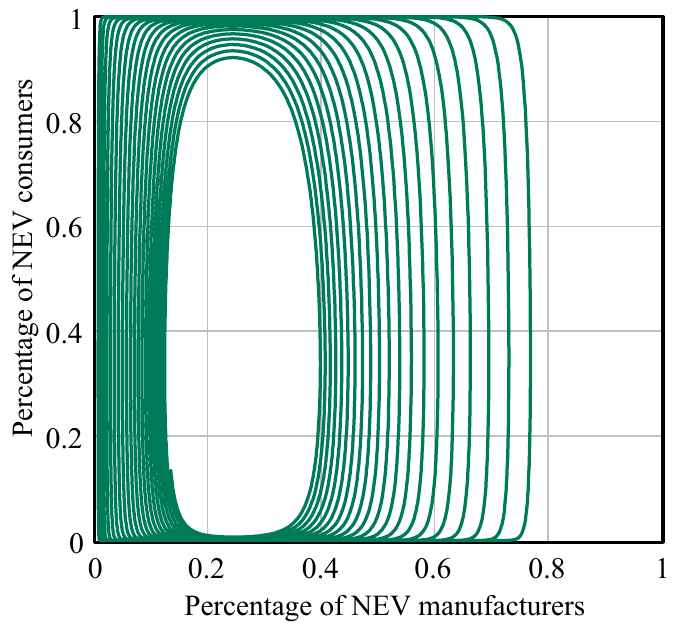"}
		\end{minipage}
	}
	\subfigure[Strategic evolution of the subject]{\label{FIG5-6b}	
		\begin{minipage}[]{0.45\linewidth}
			\centering
			\includegraphics[scale=0.4]{"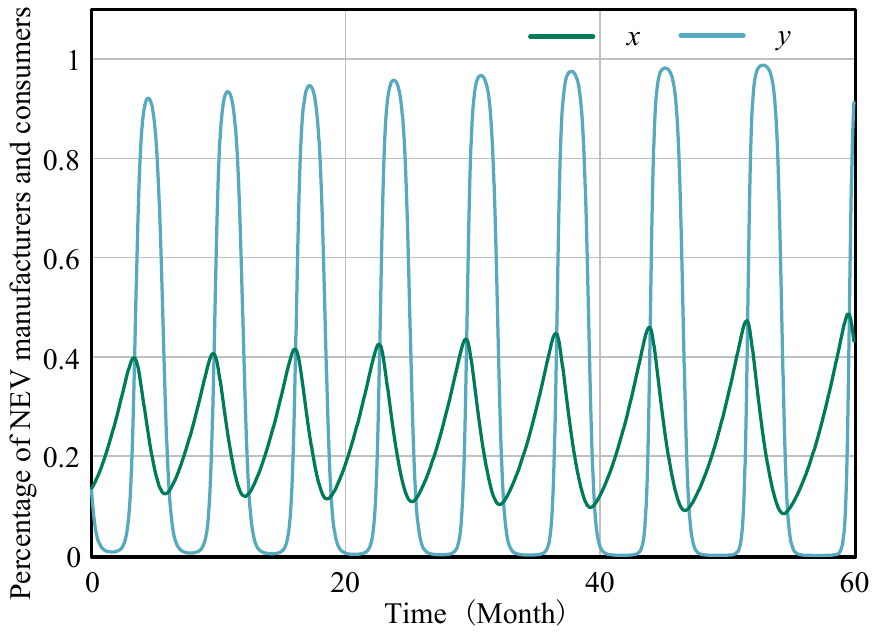"}
		\end{minipage}
	}
	\vspace{-8pt}
	\caption{\centering Model evolution trend without feedback.}
	\label{FIG5-6}
\end{figure*}

\begin{figure*}[h!]
	\subfigure[Evolutionary trajectories]{\label{FIG5-7a}	
		\begin{minipage}[]{0.45\linewidth}
			\centering
			\includegraphics[scale=0.37]{"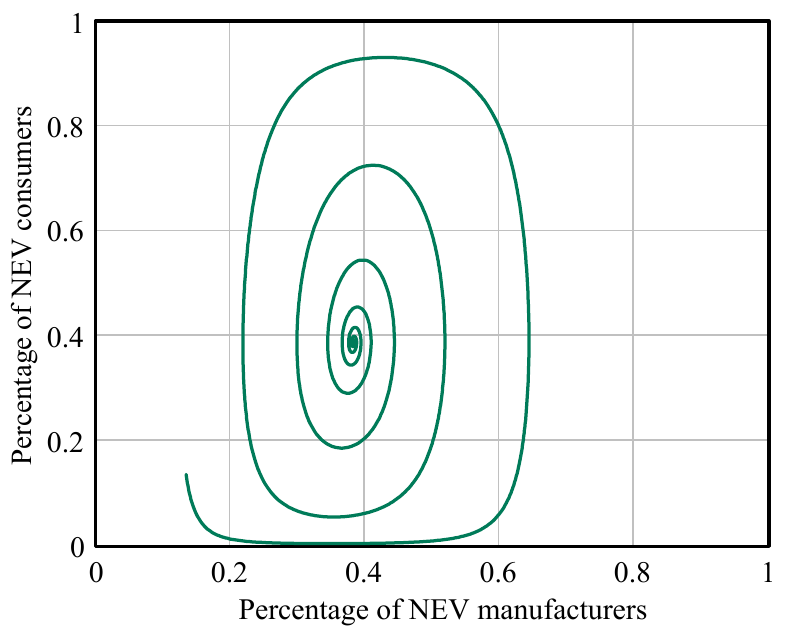"}
		\end{minipage}
	}
	\subfigure[Strategic evolution of the subject]{\label{FIG5-7b}
		\begin{minipage}[]{0.45\linewidth}
			\centering
			\includegraphics[scale=0.38]{"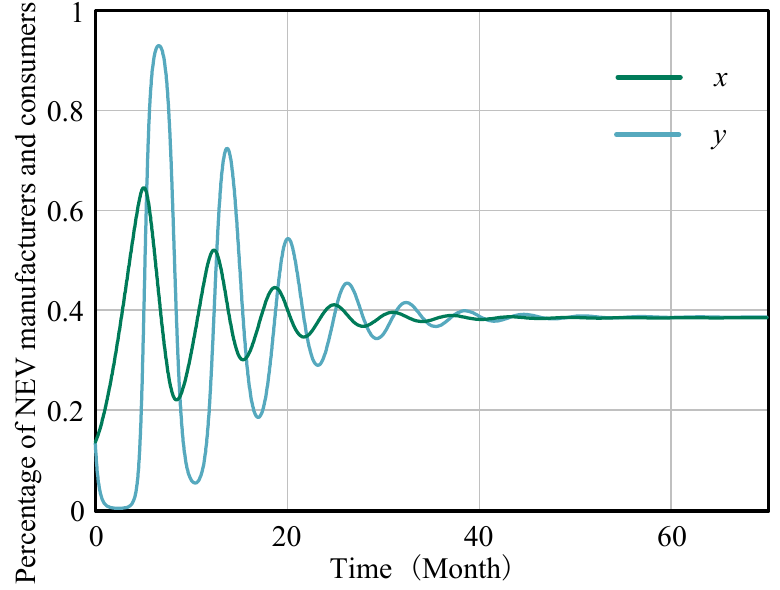"}
		\end{minipage}
	}
	\vspace{-8pt}
	\caption{\centering Model evolution trend with feedback.}
	\label{FIG5-7}
\end{figure*} 

\begin{figure*}[h!]
	\subfigure[$\sigma = 0.2$]{\label{FIG5-8a}	
		\begin{minipage}[]{0.45\linewidth}
			\centering
			\includegraphics[scale=0.38]{"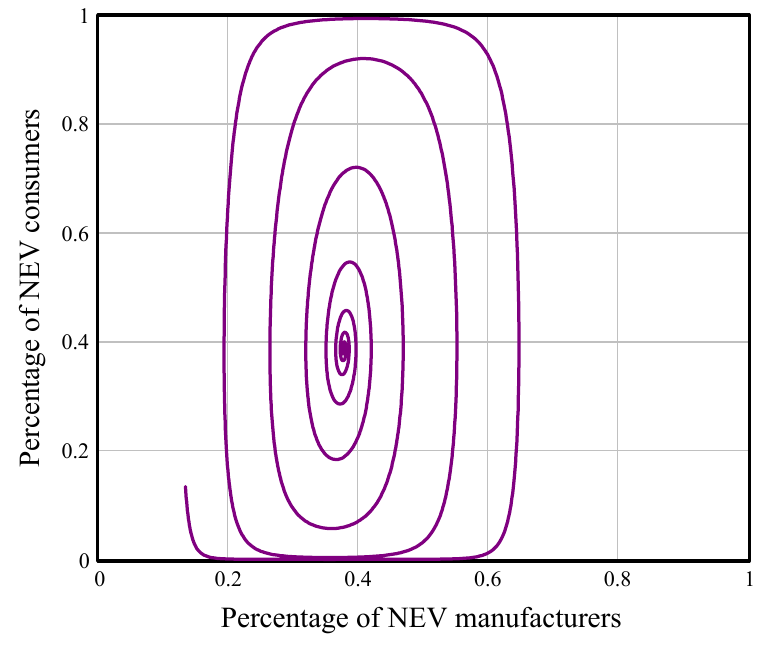"}
		\end{minipage}
	}
	\subfigure[$\sigma = 0.4$]{\label{FIG5-8b}
		\begin{minipage}[]{0.45\linewidth}
			\centering
			\includegraphics[scale=0.38]{"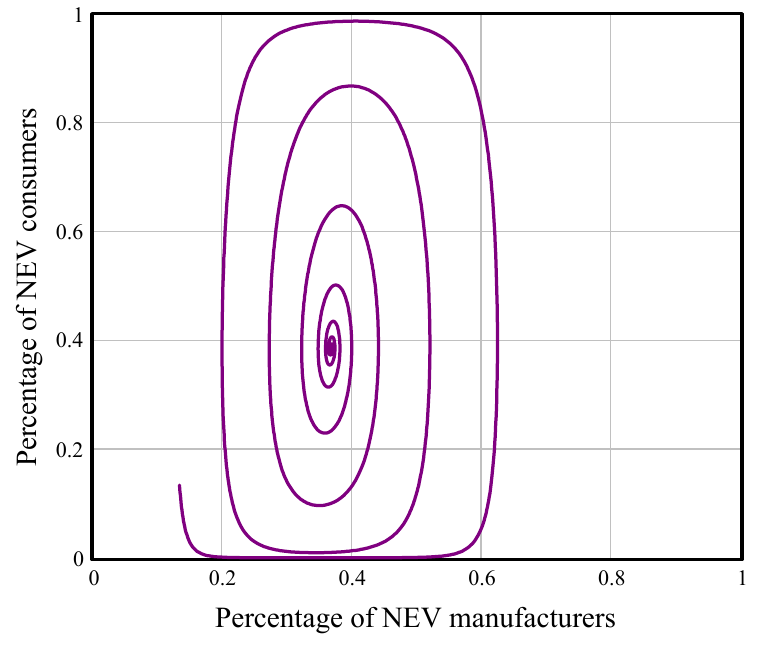"}
		\end{minipage}
	}	\subfigure[$\sigma = 0.6$]{\label{FIG5-9c}	
		\begin{minipage}[]{0.45\linewidth}
			\centering
			\includegraphics[scale=0.38]{"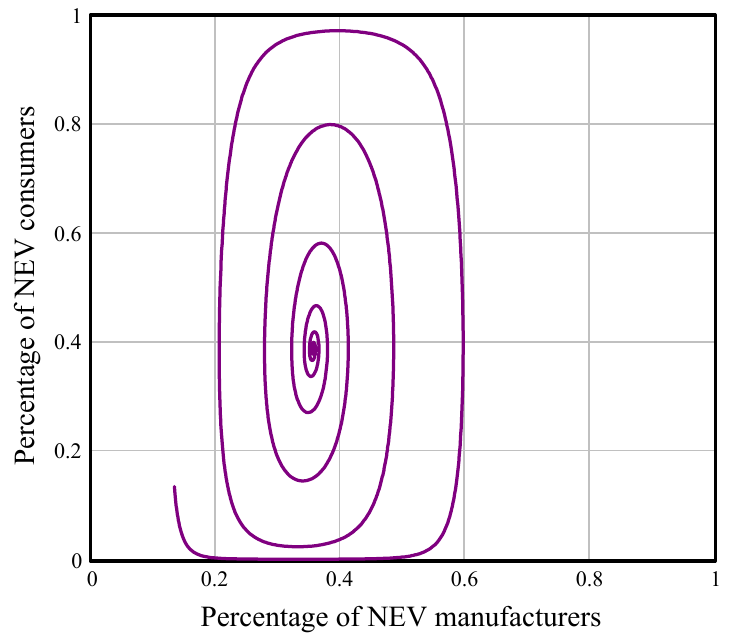"}
		\end{minipage}
	}
	\subfigure[$\sigma = 0.8$]{\label{FIG5-8d}
		\begin{minipage}[]{0.45\linewidth}
			\centering
			\includegraphics[scale=0.38]{"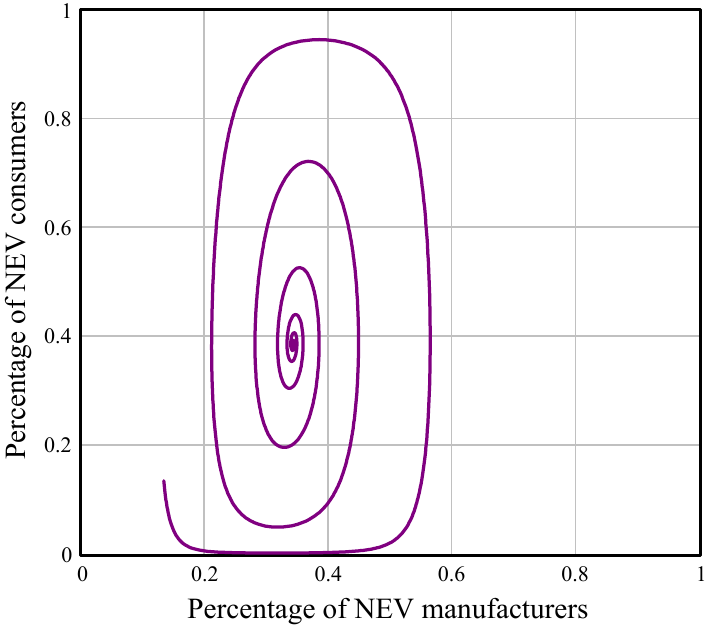"}
		\end{minipage}
	}
	\vspace{-8pt}
	\caption{\centering Model evolution trend.}
	\label{FIG5-8}
\end{figure*}

As shown in Fig.\ref{FIG5-6}, the evolutionary game model of both vehicle manufacturers and consumers without feedback has no ESS points. The strategies chosen by the vehicle manufacturer and the consumer have been fluctuating as shown in Fig.\ref{FIG5-6b} without reaching a steady state.

When the study introduced the feedback of post-purchase behaviour, the evolutionary game model gradually appeared ESS point over time, and the strategies of both subjects gradually converged to reach the Nash equilibrium stability, which is revealed in Fig.\ref{FIG5-7}. This implies that the introduction of the feedback regulation is conducive to the development of the NEVs.

The influence of the feedback regulator $\sigma$ of post-purchase behaviour is shown in Fig.\ref{FIG5-8}, where the model evolves the fastest when $\sigma=0.8$. On the contrary, the model evolves the slowest at $\sigma=0.2$. Accordingly, it is necessary for automobile manufacturers to pay attention to the after-sales service and experience of NEVs in order to reduce the depreciation of expectations brought about by consumers' post-purchase behaviour and to avoid the negative impacts caused by consumers after purchasing NEVs.

\subsection{The impact of Min-Max normalisation on evolution}

\begin{figure*}[h!]
	\subfigure{	
		\begin{minipage}[]{1\linewidth}
			\centering
			\includegraphics[scale=0.5]{"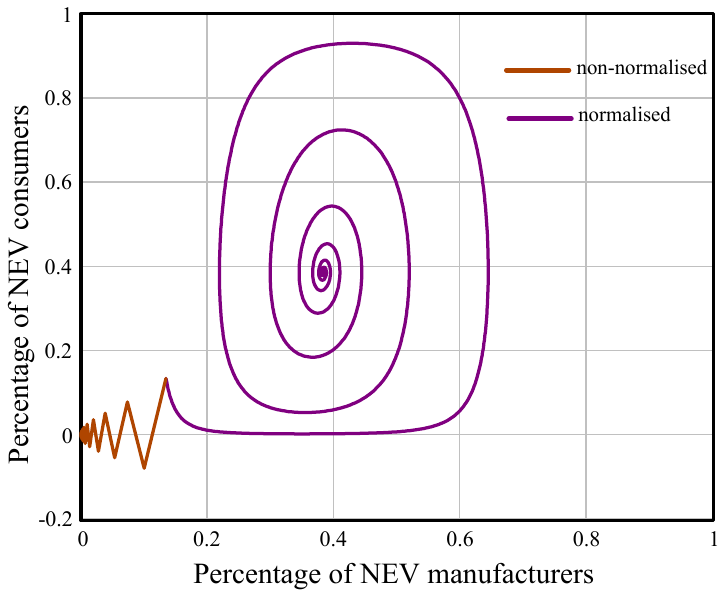"}
		\end{minipage}
	}
	\vspace{-8pt}
	\caption{\centering Model evolution trend.}
	\label{FIG5-9}
\end{figure*}

Figure \ref{FIG5-9} showed a comparison of the trajectories of the evolutionary game model before and after normalisation. Although the system has a steady state in both cases, before the normalisation of the key parameters, the system reaches ESS after only a short period of time with oscillation. And the result that both the vehicle manufacturer and the consumer choose exclusively the TFVs is to some extent not in line with reality. By comparison, the normalised pre-processed model has a relatively smooth evolutionary trajectory and the convergence endpoint is closer to reality. This suggested that the use of Min-Max normalised preprocessing data could be beneficial to improve the accuracy of the model.

\subsection{The impact of demand external stimuli on evolution}

\begin{figure*}[h!]
	\subfigure{	
		\begin{minipage}[]{1\linewidth}
			\centering
			\includegraphics[scale=0.5]{"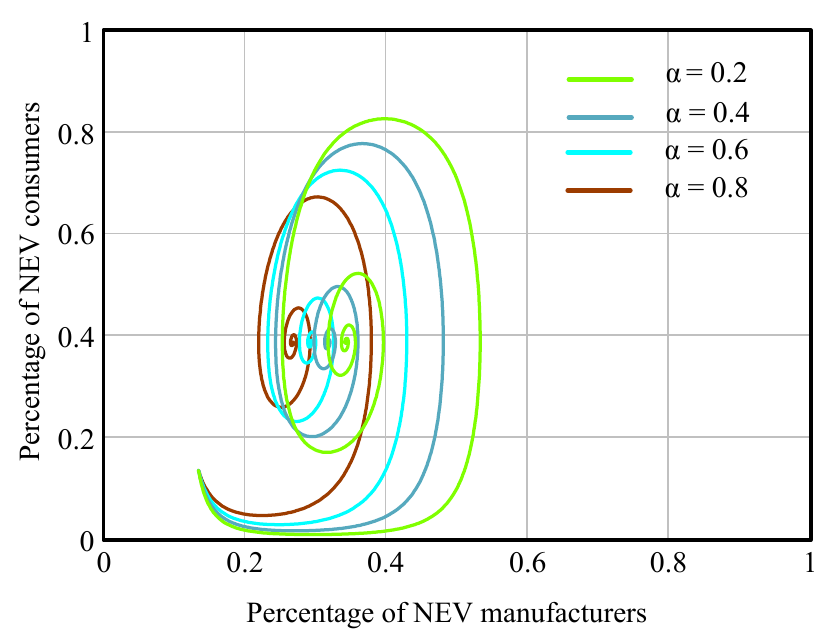"}
		\end{minipage}
	}
	\vspace{-8pt}
	\caption{\centering Model evolution trend. }
	\label{FIG5-10}
\end{figure*}

As can be seen from Fig.\ref{FIG5-10}, with the gradual increase of $a$, the evolution time of the model gradually decreases, and the system reaches the stabilisation point more quickly. However, it is worth noting that the increase of $a$ has almost no effect on the consumers who buy NEVs, but rather reduces the number of manufacturers producing NEVs, which may be due to the fact that for professional grade products such as NEVs, which are usually sold at a higher price, the consumers pay more attention to their own intrinsic needs and values and are less affected by external interference. This result also suggests that marketing activities such as promotions and advertisements should not be overly stimulating, as they may cause consumer boredom and hinder the promotion of NEVs.

\subsection{The impact ofbattery payout and replacement insurance on evolution}

\begin{figure*}[h!]
	\subfigure{	
		\begin{minipage}[]{1\linewidth}
			\centering
			\includegraphics[scale=0.5]{"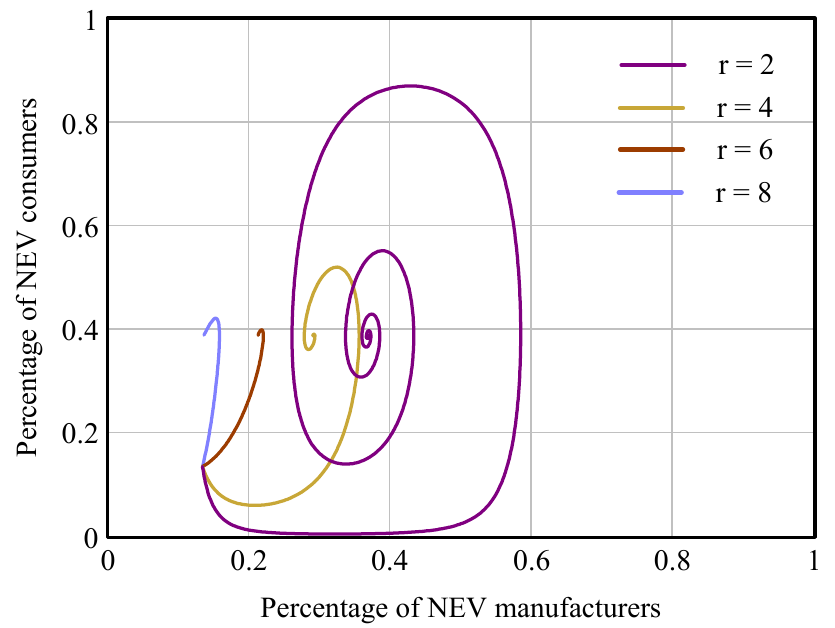"}
		\end{minipage}
	}
	\vspace{-8pt}
	\caption{\centering Model evolution trend. }
	\label{FIG5-11}
\end{figure*}

As shown in Fig.\ref{FIG5-11}, the system tends to stabilise significantly faster and the time is significantly reduced when $r$, the payout and replacement insurance for NEV batteries, gradually increases. However, when $r>6$, the evolution of the system gradually becomes slower, i.e., the time required for the model to evolve to ESS is the shortest when $r=6$. The results show that although the battery, as one of the core components of NEVs, its performance and life can affect the value assessment of NEVs by consumers, the implementation of the relevant insurance policy can not be increased in a single way, but should be adjusted according to the number of NEVs produced by the vehicle manufacturers themselves, and too high or too low settings will hinder the development of NEVs.

\begin{figure*}[h!]
	\subfigure{\label{FIG12-7a}	
		\begin{minipage}[]{1\linewidth}
			\centering
			\includegraphics[scale=0.6]{"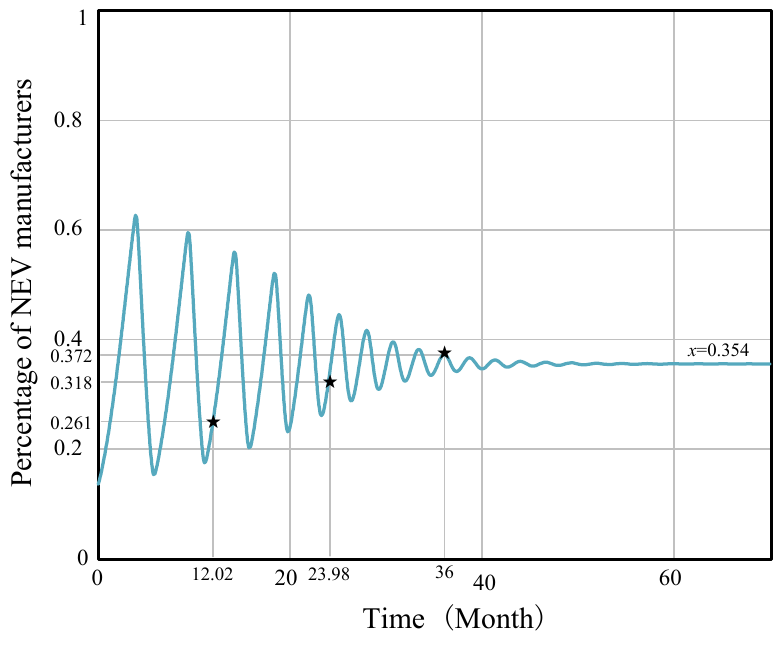"}
		\end{minipage}
	}
	\vspace{-8pt}
	\caption{\centering The trajectory of NEV development and the forecast realisation.}
	\label{FIG5-12}
\end{figure*} 

\subsection{Forecast for the promotion of NEVs}

This study not only achieved significant results in the accuracy of the evolutionary game model, but also made certain predictions about the promotion trend of NEVs at the end. Specifically, the proportion of China's NEV production in 2022 is 26.1\%\cite{2022qiche}, and the production in 2023 is 31.8\%\cite{2023qiche}. As shown in Fig.\ref{FIG5-12}, when the initial values of $x$ and $y$ are set to 0.135 and 0.134 respectively according to the development in 2021, the proportion of manufacturers producing NEVs corresponding to the horizontal axis time of 12.02 and 23.98 is 0.261 and 0.318 respectively, which fits the historical development path of NEVs in China so far. At the same time, according to the model in the time of 36 months of the vertical coordinate value of 0.372, the study predicts that China's production of NEVs manufacturers will reach 37.2\% in 2024 (that is, China's production of NEVs in 2024 can account for 37.2\%). If the current promotional measures remain unchanged, this share will eventually stabilise at 35.4\% after many years.

\begin{figure*}[h!]
	
	\subfigure{\label{FIG12-7b}
		\begin{minipage}[]{1\linewidth}
			\centering
			\includegraphics[scale=0.6]{"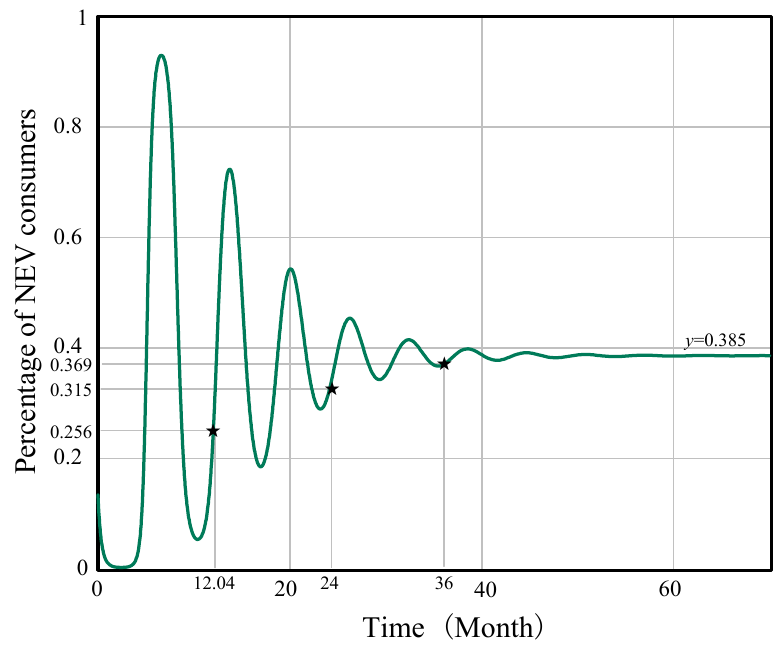"}
		\end{minipage}
	}
	\vspace{-8pt}
	\caption{\centering The trajectory of NEV development and the forecast realisation.}
	\label{FIG5-13}
\end{figure*}

Similarly, for the consumer level, the share of China's NEV sales will be 25.6\% in 2022\cite{2022qiche} and 31.5\% in 2023\cite{2023qiche}. This is consistent with the coordinate values (12.04, 0.256) and (24, 0.315) in Fig.\ref{FIG5-13}. In addition, as shown in Fig.\ref{FIG5-13}, the (36, 0.369) predicts that the sales of NEVs will reach 36.9\% in 2024, and that the sales will converge to 38.5\% after many years when the existing related promotion strategies remain unchanged.

\section{Conclusions }
\label{section6}

This study has focused on the development of NEVs, and explored feasible incentive mechanisms for the promotion of NEVs by combining theoretical knowledge from psychology and artificial intelligence to construct a more accurate evolutionary game model. The conclusions of the study are as follows:

At the theoretical level, the different mathematical approaches in the study achieved some success. Specifically, it was found that the consumer decision-making process with feedback regulation can enable the system in having a stabilisation point, and the introduction of feedback can significantly improve the stability of the inter-subjective strategy game. In addition, the Min-Max normalisation also improved the accuracy of the model, which means that it can be effective in dealing with the problem of different magnitudes in the quantitative analysis of evolutionary game models. While traditional game theories are usually based on specific forms and scales of payoff functions, these theoretical approaches can provide new ideas for future research to deal with different problems in mathematical models, which is conducive to the expansion of the application areas of game models. 

With regard to specific incentives, instead of increasing the number of consumers purchasing NEVs, external demand stimuli reduced the number of manufacturers producing NEVs. This means that vehicle manufacturers should not invest too much in promotions, advertisements and other marketing activities, as consumers who are more concerned about the actual experience may become bored. In addition, consumers' concern about damage and longevity of NEV batteries is a significant barrier to purchase. To address this issue, vehicle manufacturers can provide insurance to cover battery damage and replacement, rather than blindly launching promotional campaigns and placing excessive commercial advertisements to attract consumers. Such an insurance policy can alleviate consumer risks and concerns and increase their confidence in NEVs, thus facilitating the purchase decision. If the current promotion strategy remains unchanged, the study predicts that the production and sales ratios of NEVs in China will reach 37.2\% and 36.9\%, respectively, by 2024.

Although this study achieved the reproduction and prediction of the development path of NEVs by gradually improving the evolutionary game model, there are some limitations, which can improve some ideas and directions for future research. Firstly, as an interdisciplinary evolutionary game theory study, this study applied the remaining mathematical methods such as the consumer decision process in the modelling process, and the impact of these methods on the model was not explored. The improvement of the model accuracy by Min-Max normalisation was also at the level of numerical simulation. It may be possible to follow up with theoretical proof of the specific impact of these tools on evolutionary game models. In addition, more research can be carried out on the expectation-supply-demand game proposed in this paper, which is sufficiently feasible to be used as a base game model for exploring the emergence of cooperative behaviours in scale-free networks \cite{1xxxxxxx22} and small-world networks \cite{aaaaaa}. Finally, a wider range of incentives is also worth being discussed, such as a dynamic vehicle purchase tax or a NEV grid.

%% The Appendices part is started with the command \appendix;
%% appendix sections are then done as normal sections
%% \appendix

%% \section{}
%% \label{}

%% If you have bibdatabase file and want bibtex to generate the
%% bibitems, please use
%%

%% else use the following coding to input the bibitems directly in the
%% TeX file.

%\begin{thebibliography}{00}

%% \bibitem{label}
%% Text of bibliographic item

%\end{thebibliography}

\section*{Declarations of interest}

The authors declare that they have no known competing financial interests or personal relationships that could have appeared to influence the work reported in this paper.

\section*{CRediT authorship contribution statement}

\textbf{Tao Jin:} Conceptualization, Methodology, Software, Data Curation, Writing – Original Draft, Writing – Review \& Editing. \textbf{Yulian Jiang:} Supervision, Validation, Writing – Review \& Editing. \textbf{Xingwen Liu:} Supervision, Validation.

\section*{Acknowledgments}

This research is supported by National Nature Science Foundation (62073270), State Ethnic Affairs Commission Innovation Research Team, 2021 the Fundamental Research Funds for the Central Universities-provincial and ministerial platform construction special projects (2021NYYXS115), Innovative Research Team of the Education Department of Sichuan Province (15TD0050), and the Fundamental Research Funds for the Central Universities, Southwest Minzu University (2022NYXXS108).

%Bibliography
\bibliographystyle{unsrt}  
  
\bibliography{bibtex}

\end{document}